\documentclass[12pt]{article}
 \usepackage{graphicx}
 \usepackage{url}
\usepackage{amssymb,amsfonts,amsthm,amsmath}

\begin{document}

\begin{center}

{\bf Deducing properties of ODEs from their discretization}

\bigskip

 G.R.W. Quispel \\
Department of Mathematics, La Trobe University, VIC 3083, Australia \\
\hspace{2.0cm} r.quispel@latrobe.edu.au.  \newline  \\
 
D.I. McLaren \\
Department of Mathematics, La Trobe University, VIC 3083, Australia \\
\hspace{2.0cm} d.mclaren@latrobe.edu.au  \newline \\
 
C. Evripidou \\
Department of Mathematics, University of Cyprus \\
pamb0sd16@gmail.com

\end{center}

\section{Introduction}
\medskip \noindent Consider a polynomial ordinary differential equation (ODE)
\begin{equation}
\frac{dx}{dt} = f(x)   
\end{equation}
In this paper we address the question

\medskip \noindent Q1: What \underline{rational} conserved integral(s) and (inverse)  \underline{polynomial} volume form (if any) does the ODE possess?

\medskip \noindent Since finding rational integrals generally requires solving a nonlinear problem we propose a three step program, that, using a certain ansatz, only requires the solution of linear problems:

\medskip \noindent Step 1: Discretise the ODE using a ``suitable" method. In this paper we will use Kahan's method (but much of the following also holds for certain other birational integration methods given in the references). Compute the Jacobian determinant $J$ of the discretisation, and factorise $J$.

\medskip \noindent Step 2: Use the factors of $J$ as candidate discrete cofactors for finding discrete Darboux polynomials (DPs).

\medskip \noindent Step 3: Take the continuum limits of the discrete cofactors and DPs found in step 2. If possible, use these DPs as building  blocks for time-dependent/time-independent first integrals and preserved measure of the ODE, If one is very lucky, it may even be possible to use them to derive the exact solution to the initial value problem fore the ODE.

\section{Step 1}
Our ongoing example in this paper will be the ODE
\begin{eqnarray}
\dot{x} &=& 2-2x+xz \nonumber \\
\dot{y} &=& -y+yz \\
\dot{z} &=& -y-3z+z^2 \nonumber
\end{eqnarray}
For any quadratic ODE:
\begin{equation}
 \frac{dx_i}{dt} = \sum_{j,k}^{}a_{ijk}x_jx_k + \sum_{j}^{}b_{ij}x_j + c_i 
\end{equation}
Kahan's ``unconventional" method is defined by
\begin{equation}\label{KM}
\frac{x'_i - x_i}{h} = \sum_{j,k}^{}a_{ijk} \frac{x'_jx_k + x_jx'_k}{2} + \sum_{j}^{}b_{ij} \frac{x_j + x'_j}{2} + c_i
\end{equation}
here $x_i:=x_i(nh)$, $x'_i:=x_i((n+1)h)$, and $h$ is the timestep.

\medskip \noindent It is not hard to show that eq(\ref{KM}) can be rearranged as follows:
\begin{equation}
 \frac{x'-x}{h} = \left(I - \frac{h}{2}f'(x) \right)^{-1}f(x),   
\end{equation}
This defines the Kahan map $x_{n+1} = \phi(x_n)$ \cite{CMOQ}.

\medskip \noindent Next we compute the Jacobian determinant $J$ of $\phi$:
\begin{equation}
 J(x) = \left| \frac{\partial \phi_i(x)}{\partial x_j} \right| ,
\end{equation}
and use an algebraic manipulation package to factorise $J$.

\medskip \noindent For our example
\begin{equation}
 J(x) = \frac{K_1K_2K_3K_4}{D_1D_2^4}  ,
\end{equation}
where $K_i,D_j (i=1,\dots,4; j=1,2)$ are given in the appendix.

\section{Step 2}
\medskip \noindent Given a map $x_{n+1} = \phi(x_n)$, a polynomial $P(x)$ is a (discrete) DP of $\phi$ if there exists a (non-tautological) rational function $C$ s.t.
 \begin{equation}
P(x_{n+1}) = C(x_n)P(x_n),  
\end{equation}
where $P(x_{n+1}) = P(\phi(x_n))$ and $C$ is called the (discrete) cofactor of $P$ \cite{CEMOQTV,CEMOQT}.

\noindent\fbox{\begin{minipage}[t][1.5\height][c]{\dimexpr\textwidth-1.5\fboxsep-1.5\fboxrule\relax}
Ansatz: Given a rational map $\phi$ with Jacobian determinant
\begin{equation}
J({\bf x}) = \frac{\prod_{i=1}^{l} K_i^{b_i}({\bf x})}{ \prod_{j=1}^{m}D_j^{m_j}({\bf x})},  
\end{equation}
we try all cofactors (up to a certain polynomial degree $d$) of the form
\begin{equation}
C({\bf x}) = \pm  \frac{\prod_{i=1}^{l} K_i^{f_i}({\bf x})}{ \prod_{j=1}^{m}D_j^{g_j}({\bf x})},  
\end{equation}
where $f_i, g_j \in \mathbb{N}_0$.
\end{minipage}}

\bigskip \noindent NOTE:
\begin{enumerate}
\item There is a finite number of these co-factors up to a certain degree. For each of this finite number of co-factors, we only need to solve a linear problem (up to a chosen degree)!
\item If $C({\bf x})=J({\bf x})$, the corresponding Darboux polynomials are inverse densities of preserved measures.
\end{enumerate}
\noindent The discrete cofactors $C_i$ and corresponding DPs $P_i$ for our example are given in the first two columns of Table 1:
  \begin{table}
  \begin{center}
   \begin{tabular}{| c | c | c | c |}
     \hline
     i & $C_i$ & $P_i$ & $\bar{C}_i$ \\ \hline
    1 & $\frac{K_1}{D_2}$ & $z-y-3$ & $z$ \\ \hline
    2 & $\frac{K_2}{D_2}$ & $2z+y$ & $z-3$ \\ \hline
    3 & $\frac{K_3}{D_2}$ & $y$ & $z-1$ \\ \hline 
    4 & $\frac{K_4}{D_1D_2}$ & $x+y+z-1$ & $z-2$ \\ \hline
   \end{tabular}
  \end{center}
  \caption{}
 \end{table}

\section{Step 3}
\medskip \noindent The continuum limits $\bar{P}_i$, $\bar{C}_i$ are given by $lim_{h\rightarrow 0} P_i$ resp. $lim_{h\rightarrow 0} \frac{C_i-1}{h}$, and satisfy the ODEs
\begin{equation}
\dot{\bar{P}}_i = \bar{C_i} \bar{P_i}
\end{equation}
A useful property of the cofactor $\bar{C}_i$ is \cite{Gor}
\begin{equation}
\dot{\bar{P}}_i = \bar{C_i} \bar{P_i} \rightarrow \dot{\mathcal{P}} = \mathcal{C}\mathcal{P}
\end{equation}
where
\begin{equation}
\mathcal{P} := \prod_{i}^{} \bar{P}_i^{\alpha_i} , \mathcal{C} := \sum_{i}^{} \alpha_i \bar{C}_i
\end{equation}
This has the following implications:
\begin{enumerate}
\item $\mathcal{C}(x) = 0 \rightarrow \dot{ \mathcal{P}} = 0 \rightarrow  \mathcal{P}$ \mbox{ is a first integral}
\item $\mathcal{C}(x) = C \rightarrow \dot{ \mathcal{P}} = C \mathcal{P} \rightarrow \mathcal{P}(x(t)) = \mathcal{P}(x(0)) e^{Ct}$
\item $\mathcal{C}(x) = div f(x) \rightarrow \dot{ \mathcal{P}} = \mathcal{C} \mathcal{P} \rightarrow f$ \mbox{ preserves the measure } $\frac{dx}{\mathcal{P}(x)}$
\end{enumerate}
For our problem, the $\bar{C}_i$ are given in the last column of Table 1. (For affine DPs, $\bar{P}_i = P_i$.  For two theorems regarding affine DPs, cf \cite{CEMOQTV}).

\medskip \noindent Note that
\begin{equation}
 \bar{C_1} -  \bar{C_2} = 3,   \bar{C_1} -  \bar{C_3} = 1,   \bar{C_1} -  \bar{C_4} = 2
\end{equation}
Hence
\begin{eqnarray}
\frac{P_1}{P_2} &=& \frac{z-y-3}{2z+y} = k_2e^{3t} \\
\frac{P_1}{P_3} &=& \frac{z-y-3}{y} = k_3e^{t} \\
\frac{P_1}{P_4} &=& \frac{z-y-3}{x+y+z-1} = k_4e^{2t} 
\end{eqnarray}
and this yields 2 time-independent first integrals:
\begin{eqnarray}
I_1 &=& \frac{P_3^2}{P_1P_4} = \frac{y^2}{(z-y-3)(x+y+z-1)} \\
I_2 &=& \frac{P_3P_4}{P_1P_2} = \frac{y(x+y+z-1)}{(z-y-3)(2z+y)} 
\end{eqnarray}
Hence $f$ is integrable. Moreover $J=C_1C_2C_3C_4$ implies that $f$ preserves the measure
\begin{equation}
\frac{dxdydz}{P_1P_2P_3P_4} = \frac{dxdydz}{y(2z+y)(z-y-3)(x+y+z-1)} 
\end{equation}

\medskip \noindent Equations (15), (16), \& (17) can be combined to give the explicit solution of ODE (2)
\begin{eqnarray}
x &=& \frac{6 e^{2 t} k_{4} k_{2}+2 \left(-3 e^{t} k_{2}+\left(1+k_{2} e^{3 t}\right) k_{4}\right) k_{3}}{k_{4} \left(2 e^{3 t} k_{2} k_{3}+3 e^{2 t} k_{2}-k_{3}\right)} \\
y &=& \frac{6 k_{2} e^{2 t}}{-2e^{3t} k_2k_3 - 3e^{2t} k_2 + k_3} \\
z &=& \frac{-3 k_{3}+3 e^{2 t} k_{2}}{2 e^{3 t} k_{2} k_{3}+3 e^{2 t} k_{2}-k_{3}}
\end{eqnarray}

\medskip \noindent Acknowledgement: We are grateful to our collaborators E Celledoni, G Gubbiotti, R McLachlan, B Owren and B Tapley for many illuminating discussions.

\bigskip \noindent APPENDIX

\medskip \noindent The explicit factors of J(x) in eq (7) are:

\medskip \noindent 
\begin{eqnarray*}
K_1 &=&  -1/4\,{h}^{2}x_{{2}}+3/4\,{h}^{2}x_{{3}}-3/4\,{h}^{2}-1/2\,hx_{{3}}-h+1 \\
K_2 &=& -1/4\,{h}^{2}x_{{2}}-1/4\,{h}^{2}x_{{3}}-3/4\,{h}^{2}-1/2\,hx_{{3}}+h+1 \\
K_3 &=& -1/4\,{h}^{2}x_{{2}}-3/4\,{h}^{2}x_{{3}}+3/4\,{h}^{2}-1/2\,hx_{{3}}+2\,h+1 \\
K_4 &=& 1/8\,{h}^{3}x_{{2}}x_{{3}}-1/8\,{h}^{3}{x_{{3}}}^{2}-1/4\,{h}^{3}x_{{2
}}+{\frac {7\,{h}^{3}x_{{3}}}{8}}+1/4\,{h}^{2}{x_{{3}}}^{2}-3/4\,{h}^{
3} \\ & &-1/4\,{h}^{2}x_{{2}}-1/4\,{h}^{2}x_{{3}}-5/4\,{h}^{2}-hx_{{3}}+h+1 \\
D_1 &=& -1/2\,hx_{{3}}+h+1 \\
D_2 &=& 1/2\,{h}^{2}{x_{{3}}}^{2}+1/4\,{h}^{2}x_{{2}}-5/4\,{h}^{2}x_{{3}}+3/4
\,{h}^{2}-3/2\,hx_{{3}}+2\,h+1
\end{eqnarray*}


\end{document}